\documentclass{article}
\usepackage{latexsym}
\usepackage{amssymb,amsmath}
\usepackage{amscd}

\newtheorem{theorem}{Theorem}[section]
\newtheorem{proposition}[theorem]{Proposition}
\newtheorem{lemma}[theorem]{Lemma}
\newtheorem{corollary}[theorem]{Corollary}

\newtheorem{remark}[theorem]{Remark}

\newcommand{\be}{\begin{eqnarray}}
\newcommand{\ee}{\end{eqnarray}}
\newcommand{\ben}{\begin{eqnarray*}}
\newcommand{\een}{\end{eqnarray*}}
\numberwithin{equation}{section}
\newcommand{\NN}{\mathbb N}
\newcommand{\ZZ}{\mathbb Z}

\newcommand{\FF}{\mathbb F}
\newcommand{\PP}{\mathbb P}

\def\cO{\mathcal O}

\def\cC{\mathcal C}
\def\cD{\mathcal D}

\def\cH{\mathcal H}

\def\cL{\mathcal L}

\def\cO{\mathcal O}

\def\cX{\mathcal X}

\def\l{\ell}

\def\K{\mathbb{K}}

\def\min{{\rm min}}

\def\fqq{{\mathbb F}_{q^2}}

\def\e{{\epsilon}}
\begin{document}
\title{AG Codes on certain Maximal Curves}

\author{Stefania~Fanali~and~Massimo~Giulietti
\thanks{S. Fanali and M. Giulietti are with the Dipartimento di Matematica e Informatica, Universit\`a di  Perugia,
Via Vanvitelli 1, 06123, Perugia,
Italy (e-mail: stefania.fanali@dipmat.unipg.it;  giuliet@dipmat.unipg.it)}
\thanks{This research was performed within the activity of GNSAGA of the
Italian INDAM.}}
\maketitle

\begin{abstract}
Algebraic Geometric codes associated to a recently discovered class of maximal curves are investigated.  As a result, some linear codes with better parameters with respect to the previously known ones are discovered, and $70$ improvements on MinT's tables \cite{MINT} are obtained.
\end{abstract}


\section{Introduction}
Algebraic Geometric codes  (AG codes) are linear error-correcting
codes constructed from algebraic curves \cite{GO1,GO2}. Roughly
speaking, the parameters of an AG code are good when the
underlying curve has many rational points with respect to its
genus. AG codes from specific curves with many points, such as the
Hermitian curve and its quotients, the Suzuki curve, and the Klein
quartic, have been the object of several works, see e.g.
\cite{HA,HS,M0,M1, ST, TI,XC,XL,YK} and the references therein.

In this paper we provide an explicit construction of one-point AG
codes from the GK curves, together with some results on the
permutation automorphism groups of such codes. The GK curves are
defined over any finite field of order $q^2$ with $q={\bar q}^3$, and they
are maximal curves in the sense that the number of their
$\fqq$-rational points attains the  Hasse-Weil upper bound
$$q^2+1+2gq$$
where $g$ is the genus of the curve. Significantly, for $q>8$, GK
curves are the first known examples of  maximal curves which are
proven not to be $\FF_{q^2}$-covered by the Hermitian curve, see
\cite{GK}.

Interestingly, some of the codes constructed in this paper have
better parameters compared with the known linear error-correcting
codes, see Section \ref{sec5}. More precisely, we obtain an improvement on
the best known mininum distance for linear codes over the finite
field with $64$ elements in the following cases
$$
\begin{array}{cc}
{\rm length} & {\rm codimension} \\
200-224 & 20 \\
210-224 & 22 \\
210-224 & 23 \\
210-224 & 28 \\
\end{array}
$$
(see Theorem \ref{new}).

The paper is organized as follows. In Section \ref{sec2}, some basic facts on
AG codes and maximal curves are recalled.
In Section \ref{sec3}, we introduce the
GK curves, by summarizing some of the results in \cite{GK}. The
Weierstrass semigroup at rational points of GK curves is
investigated in Section \ref{Weie}.  Finally, certain AG codes associated to
the GK curves are constructed and their parameters are discussed for
${\bar q}=2,3$, see Section \ref{sec5}.

\section{Preliminaries}\label{sec2}

\subsection{Curves} Throughout the paper, by a curve we mean a projective, geometrically irreducible, non-singular algebraic curve defined over a finite field.  Let $q$ be a prime power,  and let $\cX$ be  a curve defined over the finite field $\FF_{q^{2}}$ of order $q^{2}$. Let $g$ be the genus of $\cX$.
Henceforth, the following notation is used:
\begin{enumerate}
\item[$\bullet$]  $\cX(\FF_{q^{2}})$ (resp. $\FF_{q^{2}}(\cX)$) denotes the set of $\FF_{q^{2}}$-rational points (resp. the field of $\FF_{q^{2}}$-rational functions) of $\cX$.
\item[$\bullet$] For $f\in\FF_{q^{2}}(\cX)$,  $div(f)$ (resp. $div_{\infty}(f)$) denotes the divisor (resp. the pole divisor) of  $f$.
\item[$\bullet$] Let $P$ be a point of $\cX$. Then $v_{P}$ (resp. $H(P)$) stands for the valuation (resp. for the Weierstrass non-gap semigroup) associated to $P$. The $i$th non-gap at $P$ is denoted as $m_i(P)$.
\item[$\bullet$] Let $D$ be a divisor on $\cX$ and $P\in\cX$. Then $deg(D)$ denotes the degree of $D$,  $supp(D)$ the support of $D$, and  $v_{P}(D)$ the coefficient of $P$ in $D$. For $D$ an $\FF_{q^{2}}$-divisor, let
$$L(D):=\left\{f\in\FF_{q^{2}}(\cX) | div(f)+D\geq 0 \right\},$$
$$l(D):=dim_{\FF_{q^{2}}}(L(D)).$$
\item[$\bullet$] The symbol "$\sim$" denotes linear equivalence of divisors.
\item[$\bullet$] The symbol $g_{d}^{r}$ stands for a linear series of projective dimension $r$ and degree $d$.
\end{enumerate}

\subsection{One-point AG Codes and Improved AG Codes}
Let $\cX$ be a curve, let $P_{1},P_{2},\ldots,P_{n}$ be $\FF_{q^{2}}$-rational points of $\cX$, and let $D$ be the divisor $P_{1}+P_{2}+\ldots+P_{n}$. Furthermore, let  $G$ be some other divisor that has support disjoint from $D$.
The AG code $C(D,G)$ of length $n$ over $\FF_{q^{2}}$ is the image of the linear map $\alpha : L(G)\rightarrow\FF_{q^{2}}^{n}$ defined by $\alpha(f)=(f(P_{1}),f(P_{2}),\ldots,f(P_{n}))$. If $n$ is bigger than $deg(G)$, then $\alpha$ is an embedding, and the dimension $k$ of  $C(D,G)$ is equal to $\ell(G)$. The Riemann-Roch theorem makes it possible to estimate the parameters of $C(D,G)$. In particular, if $2g-2<deg(G)<n$, then $C(D,G)$ has dimension $k=deg(G)-g+1$ and minimum distance $d\geq n-deg(G)$, see e.g. \cite[Theorem~2.65]{HLP}. A generator matrix $M$ of $C(D,G)$ is
$$M=\left(\begin{array}[pos]{ccc}
f_{1}(P_{1}) & \ldots & f_{1}(P_{n}) \\
\vdots & \ldots & \vdots \\
f_{k}(P_{1}) & \ldots & f_{k}(P_{n})\\
\end{array}\right),$$
where $f_{1},f_{2},\ldots,f_{k}$ is an $\FF_{q^{2}}$-basis of $L(G)$.
The dual code $C^{\perp}(D,G)$ of $C(D,G)$ is an AG code with dimension $n-k$ and minimum distance greater than or equal to $deg(G)-2g+2$.
When $G=\gamma P$  for an $\FF_{q^{2}}$-rational $P$ point of $\cX$, and a positive integer $\gamma$, AG codes $C(D,G)$ and $C^{\perp}(D,G)$ are referred to as one-point AG codes. We recall some results on the minimum distance of one-point AG codes. By \cite[Theorem~3]{GKL}, we can assume that $\gamma$ is a non-gap at $P$. Let
$$H(P)=\left\{\rho_{1}=0<\rho_{2}<\ldots\right\},$$
and set $\rho_{0}=0$.
Let $f_{\l}$ be a rational function such that $div_{\infty}(f_{\l})=\rho_{\l}P$, for any $\l\geq1$. Let $D=P_{1}+P_{2}+\ldots+P_{n}$. Let also
\begin{equation}\label{hl}
h_{\l}=(f_{\l}(P_{1}),f_{\l}(P_{2}),\ldots,f_{\l}(P_{n}))\in\FF_{q^{2}}^{n}.
\end{equation}
%

Set
$$\nu_{\l}:=\#\left\{(i,j)\in\NN^{2} : \rho_{i}+\rho_{j}=\rho_{\l+1}\right\}$$
for any $\l\geq0$. Denote with $C_{\l}(P)$ the dual of the AG code $C(D,G)$, where $D=P_{1}+P_{2}+\ldots+P_{n}$, and $G=\rho_{\l}P$.

\begin{lemma}
\cite[Proposition~4.11]{HLP}
If $y\in C_{\l}(P)\setminus C_{\l+1}(P)$, then the weight of $y$ is greater than or equal to $\nu_{\l}$.
\end{lemma}

The integer
$$d_{ORD}(C_{\l}(P)):=min\left\{\nu_{m} : m\geq\l\right\}$$
is called the order bound or the Feng-Rao designed minimum distance of $C_{\l}(P)$.

\begin{theorem}\label{dist}
\cite[Theorem~4.13]{HLP} The minimum distance $d(C_{\l}(P))$ of $C_{\ell}(P)$ satisfies
$$d(C_{\l}(P))\geq d_{ORD}(C_{\l}(P)).$$

\end{theorem}

\begin{theorem}
\cite[Theorem~5.24]{HLP} Let
$$c:=\textrm{max}\left\{m\in\ZZ | m-1\notin H(P)\right\}.$$
Then $d_{ORD}(C_{\l}(P))\geq\l+1-g$. If $\l\geq 2c-g-1$, then $\nu_\ell=\ell+1-g$  and hence equality $d_{ORD}(C_{\l}(P))\geq\l+1-g$ holds.
\end{theorem}

Let $d$ be an integer greater than $1$. The improved AG code $\tilde{C}_d(P)$ is the code
$$\tilde{C}_d(P):=\left\{x\in\FF_{q^{2}}^{n} : \left\langle x,h_{i+1}\right\rangle=0 \textrm{ for all } i \textrm{ such that } \nu_{i}<d\right\},$$
see \cite[Def. 4.22]{HLP}.

\begin{theorem}\cite[Proposition 4.23]{HLP}\label{dist2}
Let
$$r_{d}:=\#\left\{i\geq0 : \nu_{i}<d\right\}.$$
Then $\tilde{C}_d(P)$ is an $\left[n,k,d'\right]$-code, where
$k\geq n-r_{d}$, and $d'\geq d$.
\end{theorem}

\subsection{Weierstras point theory \cite{SV}}
Let $\cD$ be a $g_{d}^{r}$ base-point-free $\FF_{q^{2}}$-linear series on a curve $\cX$.
For a point $P\in \cX$, let
$$
j_{0}(P)=0<j_{1}(P)<\ldots<j_{r}(P)\leq d
$$
be the $(\cD,P)$-orders, that is, the integers $j$ such that there exists a divisor $D\in \cD$ with $v_P(D)=j$.
This sequence is the same for all but finitely many points. The finitely many points $P$ where exceptional $(\cD, P)$-orders occur, are called the $\cD$-Weierstrass points of $\cX$.
Let $\e_{0}<\e_{1}<\ldots<\e_{r}$ denote the sequence of the  $(\cD, Q)$-orders for a generic point $Q\in\cX$. Then $\e_{i}\leq j_{i}(P)$, for each $i=0,1,\ldots r$ and for any point $P$. The ramification divisor of $\cD$ is  a divisor $R$ whose support consists  exactly of the $\cD$-Weierstrass points, and such that
$deg(R)=(\e_{0}+\e_{1}+\ldots+\e_{r})(2g-2)+(r+1)d$.


\subsection{Maximal Curves} A curve $\cX$ is called $\FF_{q^{2}}$-maximal if the number of its $\FF_{q^{2}}$-rational points attains
the Hasse-Weil upper bound, that is,
$$\#\cX(\FF_{q^{2}})=q^{2}+1+2gq,$$
where $g$ is the genus of $\cX$.

A key tool for the investigation of maximal curves is Weiestrass Points theory.
The Frobenius linear series of a maximal curve $\cX$ is the complete linear series $\cD=|(q+1)P_0|$, where $P_0$ is any $\FF_{q^2}$-rational point of $\cX$.
The next result provides a relationship between $\cD$-orders and non-gaps at points of $\cX$.

\begin{proposition}\label{FGT4}
\cite[Proposition~1.5]{FGT}
Let $\cX$ be a maximal curve over $\FF_{q^{2}}$, and let $\cD$ be the Frobenius linear series of $\cX$. Then
\begin{enumerate}
\item[(i)] For each point $P$ on $\cX$, we have $\l(qP)=r$, i.e., 
$$0<m_{1}(P)<\ldots<m_{r-1}(P)\leq q<m_{r}(P).$$
\item[(ii)] If $P$ is not rational over $\FF_{q^{2}}$, the $\cD$-orders at the point $P$ are
$$0\leq q-m_{r-1}(P)<\ldots<q-m_{1}(P)<q.$$
\item[(iii)] If $P$ is rational over $\FF_{q^{2}}$, the $(\cD, P)$-orders are
$$0<q+1-m_{r-1}(P)<\ldots<q+1-m_{1}(P)<q+1.$$
In particular, if $j$ is a $\cD$-order at a rational point $P$, then $q+1-j$ is a non-gap at $P$.
\item[(iv)] If $P\in\cX(\FF_{q^{2}})$, then $q$ and $q+1$ are non-gaps at $P$.
\end{enumerate}
\end{proposition}

Maximal curves are characterized by the so-called Natural Embedding Theorem.

\begin{theorem}\label{teo_imm}
\cite[Theorem~10.22]{HKT}
Every $\FF_{q^{2}}$-maximal curve $\cX$ of genus $g\geq0$ is isomorphic over $\FF_{q^{2}}$ to a curve of $\PP^{m}(\bar{\FF}_{q^{2}})$ of degree $q+1$ lying on a non-degenerate Hermitian variety $\cH_{m}$ defined over $\FF_{q^{2}}$.
\end{theorem}

The dimension $m$ in Theorem \ref{teo_imm} is less than or equal to the dimension $r$ of the Frobenius linear series of $\cX$. Also,
by \cite[Theorem~10.22]{HKT}, the osculating hyperplane of $\cX$ at any point $P\in\cX$ coincides with the tangent hyperplane at $P$ to the non-degenerate Hermitian variety $\cH_{m}$ in which $\cX$ lies.

%
\section{GK-Curves}\label{sec3}
Let $q=\bar{q}^3$, where $\bar{q}\ge 2$ is a prime power.
The GK-curve over $\FF_{q^2}$ is the curve of ${\mathbb P}^3(\bar{\FF}_{q^{2}})$ with affine equations
\begin{equation}\label{equazioni}
\left\{\begin{array}{ll}
Z^{\bar{q}^{2}-\bar{q}+1}=Y h(X)\\
X^{\bar{q}}+X=Y^{\bar{q}+1}
\end{array}\right.,
\end{equation}
where $h(X)=\sum_{i=0}^{\bar{q}}(-1)^{i+1}X^{i(\bar{q}-1)}$.
We first recall some important proprieties of this curve, for which we refer to \cite[Section~2]{GK}.
The curve $\cX$ is absolutely irreducible, non-singular, and it lies on the Hermitian surface $\cH_{3}$ with affine equation
$$X^{\bar{q}^{3}}+X=Y^{\bar{q}^{3}+1}+Z^{\bar{q}^{3}+1}.$$
Hence, by \cite[Theorem~10.31]{HKT}, $\cX$ is $\FF_{q^{2}}$-maximal. Significantly, for $q>8$, $\cX$ is the only known curve that is maximal but not $\FF_{q^2}$-covered by the Hermitian curve $\cH_{2}$ defined over $\FF_{q^{2}}$ (see \cite[Theorem~5]{GK}). The genus of
$\cX$ is
$$g=\frac{1}{2}(\bar{q}^{3}+1)(\bar{q}^{2}-2)+1.$$

In order to investigate one-point AG codes associated to $\cX$, we need to describe the Weierstrass semigroup  $H(P)$ associated to an $\FF_{q^{2}}$-rational point $P\in\cX$. In the rest of this section we establish some general properties of $H(P)$. 

Let $\Lambda$ be the cyclic group consisting of all collineations
$$g_{u} : (T:X:Y:Z)\longmapsto(uT:uX:uY:Z),$$
with $u^{\bar{q}^{2}-\bar{q}+1}=1$. Clearly $\Lambda$ is a projective group preserving $\cX$. It is easily seen that a plane model for the quotient curve $\cX/\Lambda$ has equation $X^{\bar{q}}+X=Y^{\bar{q}+1}$. Consider the projection $\psi : \cX\rightarrow\cX/\Lambda$.
Let $P=(X_{P},Y_{P},Z_{P})$ be any affine $\FF_{q^{2}}$-rational point of $\cX$. Then, either $\psi$ is fully ramified at $P$, or $\psi$ splits completely at $\bar{P}=\psi(P)$, according to whether $\psi(P)$ is an $\FF_{\bar{q}^{2}}$-rational point of $\cX/\Lambda$ or not, that is, whether $Z_{P}=0$ or not.

In the former case, let $aX+bY+c=0$ be an equation of any line through $\psi(P)$, distinct from the tangent of $\cX/\Lambda$ at $\psi(P)$. Then
$$v_{P}(ax+by+c)=(\bar{q}^{2}-\bar{q}+1)v_{\psi(P)}(ax+by+c)=\bar{q}^{2}-\bar{q}+1$$
holds.
Hence, by \cite[Prop.~10.6~(IV)]{HKT}, $\bar{q}^{3}+1-(\bar{q}^{2}-\bar{q}+1)=\bar{q}^{3}-\bar{q}^{2}+\bar{q}\in H(P)$.

By \cite[Proposition 1]{GK}, together with \cite[Theorem 7]{GK}, $\bar{q}^{3}-\bar{q}^{2}+\bar{q},\bar{q}^3,\bar{q}^3+1$ are actually a set of generators for $H(P)$. The same holds when $P$ is the only infinite point of $\cX$.
\begin{proposition}\label{0}
If either $P$ is the only infinite point of $\cX$ or  $P=(X_{P}, Y_{P}, 0)\in\cX(\FF_{{q}^{2}})$, then the Weierstrass semigroup at $P$ is the subgroup generated by $\bar{q}^{3}-\bar{q}^{2}+\bar{q}$, $\bar{q}^{3}$, and $\bar{q}^{3}+1$.
\end{proposition}

Assume now that $\psi(P)$ is a non-$\FF_{\bar{q}^{2}}$-rational point of $\cX/\Lambda$. Let $aX+bY+c=0$ be an equation of the tangent line of $\cX/\Lambda$
at $\psi(P)$. The order of contact of the tangent line to $\cX/\Lambda$ at a non-$\FF_{\bar{q}^{2}}$-rational point is equal to $\bar{q}$ (see e.g. \cite[p. 302]{HKT}). Then
$$v_{P}(ax+by+c)=v_{\psi(P)}(ax+by+c)=\bar{q}$$
holds.
Again by \cite[Prop.~10.6~(IV)]{HKT}, $\bar{q}^{3}-\bar{q}+1\in H(P)$. Then the following result is obtained.

\begin{proposition}
If $P=(X_{P}, Y_{P}, Z_{P})\in\cX(\FF_{{q}^{2}})$ is such that $Z_{P}\neq0$, then the Weierstrass semigroup at $P$ contains the subgroup generated by $\bar{q}^{3}-\bar{q}+1$, $\bar{q}^{3}$, and $\bar{q}^{3}+1$.
\end{proposition}

Providing a general description of $H(P)$ for affine points $P$ of $\cX$ with $Z_P\neq 0$ seems to be quite a difficult task. In the next section, this will be done for the cases $\bar{q}=2,3$.

We end this section by describing the automorphism group $Aut(\cX)$ of $\cX$, together with its action on the set of $\FF_{q^{2}}$-rational points of $\cX$.
\begin{theorem}\cite[Theorem~6]{GK}
\label{automorf} The automorphism group of $\cX$ has order $\bar{q}^3(\bar{q}^3+1)(\bar{q}^2-1)(\bar{q}^2-\bar{q}+1)$, and has a normal subgroup isomorphic to $SU(3,\bar{q})$. If ${\rm gcd}(3,\bar{q}+1)=1$ then  $Aut(\cX)$ is isomorphic to the direct product of $SU(3,\bar{q})$ and a cyclic group of order $\bar{q}^2-\bar{q}+1$. If
${\rm gcd}(3,\bar{q}+1)=3$ then  $Aut(\cX)$ has a normal subgroup of index $3$ which is isomorphic to the direct product of $SU(3,\bar{q})$ and a cyclic group of order $(\bar{q}^2-\bar{q}+1)/3$. 
\end{theorem}

\begin{theorem}\cite[Theorem~7]{GK}
\label{action} The set of $\FF_{q^{2}}$-rational points of $\cX$ splits into two orbits under the action of $Aut(\cX)$. One orbit, say $\cO_1$, has size $\bar{q}^3+1$ and consists of the points $(X_P,Y_P,0)\in \cX(\FF_{q^{2}})$, together with the infinite point $X_\infty=(0:1:0:0)$. The other orbit, say $\cO_2$, has size
$\bar{q}^3(\bar{q}^3+1)(\bar{q}^2-1)$ and consists of the points $(X_P,Y_P,Z_P)\in \cX(\FF_{q^{2}})$ with $Z_P\neq 0$. Furthermore, $Aut(\cX)$ acts on $\cO_1$ as $PGU(3,\bar{q})$ in its doubly transitive permutation representation.
\end{theorem}

\begin{corollary}\label{stabilizer} For a point $P\in \cO_1$, the stabilizer of $P$ under the action of $Aut(\cX)$ has size $\bar{q}^3(\bar{q}^2-1)(\bar{q}^2-\bar{q}+1)$, and it acts transitively on the points of $\cO_1 \setminus \{P\}$.
For a point $P\in \cO_2$, the stabilizer of $P$ under the action of $Aut(\cX)$ has size $(\bar{q}^2-\bar{q}+1)$.
\end{corollary}

\section{The Weierstrass semigroup at an $\FF_{q^{2}}$-rational point of the GK curves}\label{Weie}
In this section we describe the Weierstrass semigroup $H(P)$ at any $\FF_{q^{2}}$-rational point $P$ of the GK curves for the cases $\bar{q}=2,3$. Also, for each non-gap $m$ we provide a rational function $f$ such that $div_\infty(f)=mP$.

First, consider the point at infinity $X_{\infty}=(0:1:0:0)$ of $\cX$. By Proposition~\ref{0},  $H(X_{\infty})=\left\langle \bar{q}^{3}-\bar{q}^{2}+\bar{q}, \bar{q}^{3}, \bar{q}^{3}+1\right\rangle$.
Taking into account that the osculating plane of $\cX$ at $X_{\infty}$ is the plane with equation $T=0$, and that $X_{\infty}\in\cX(\FF_{q^{2}})$,
\begin{equation}\label{eqqu1}
div_{\infty}(x)=(\bar{q}^{3}+1)X_{\infty}
\end{equation}
holds. Moreover, by the equations of $\cX$ it follows that
\begin{equation}\label{eqqu2}
div_{\infty}(y)=(\bar{q}^{3}-\bar{q}^{2}+\bar{q})X_{\infty},\,\,\,\,div_{\infty}(z)=\bar{q}^{3}X_{\infty}.
\end{equation}
It should be noted that taking into account \cite[Theorem 7]{GK} one can easily construct rational functions corresponding to the non-gaps $\bar{q}^{3}-\bar{q}^{2}+\bar{q}, \bar{q}^{3}, \bar{q}^{3}+1$ at any point $P=(a,b,0)\in\cX(\FF_{q^{2}})$.

Now fix $P=(a,b,c)\in\cX(\FF_{q^{2}})$, with $c\neq0$, and consider the planes
$\pi_{1} : T=0$,
$\pi_{2} : Y-bT=0$, $\pi_{3} : -a^{\bar{q}}T-X+b^{\bar{q}}Y=0$, and $\pi_{4}:
-a^{\bar{q}^{3}}T-X+b^{\bar{q}^{3}}Y+c^{\bar{q}^{3}}Z=0$. It is straightforward to
check that these planes meet $\cX$ at $P$ with multiplicity $0$,
$1$, $\bar{q}$, and $\bar{q}^{3}+1$ respectively. By Proposition~\ref{FGT4},
$\pi_{1},\ldots,\pi_{4}$ correspond to the rational functions
associated to the following non-gaps at $P$: $\bar{q}^{3}+1$, $\bar{q}^{3}$,
$\bar{q}^{3}-\bar{q}+1$, and $0$. Let $\phi$ be the linear transformation
$$\phi(T:X:Y:Z)=$$
$$(-a^{\bar{q}^{3}}T-X+b^{\bar{q}^{3}}Y+c^{\bar{q}^{3}}Z: T: -bT+Y: -a^{\bar{q}}T-X+b^{\bar{q}}Y).$$
Note that the planes $\phi(\pi_{1}),\ldots,\phi(\pi_{4})$ are the planes $X=0$, $Y=0$, $Z=0$, and $T=0$. Also, $\phi(P)=X_{\infty}$. The equations of $\phi(\cX)$ are
$$\displaystyle{\left(\frac{T}{c^{\bar{q}^{3}}}+cX+\left(\frac{b^{\bar{q}}-b^{\bar{q}^{3}}}{c^{\bar{q}^{3}}}\right)Y-\frac{Z}{c^{\bar{q}^{3}}}\right)^{\bar{q}^{2}-\bar{q}+1} =}$$ 
$$\left(bX+Y\right) h(aX+b^{\bar{q}}Y-Z)$$ 
and 
$$X\left(aX+b^{\bar{q}}Y-Z\right)^{\bar{q}}+ X^{\bar{q}}\left(aX+b^{\bar{q}}Y-Z\right)
=\left(bX+Y\right)^{\bar{q}+1}.$$

The rational functions $x$, $y$, $z$ correspond to the first non-gaps at $X_{\infty}$, namely
\begin{equation}\label{ng09}
v_{X_{\infty}}(x)=-(\bar{q}^{3}+1),v_{X_{\infty}}(y)=-\bar{q}^{3},v_{X_{\infty}}(z)=-(\bar{q}^{3}-\bar{q}+1).
\end{equation}
We now look for non-gaps at $X_{\infty}$ which are not in the subgroup generated by $\bar{q}^{3}-\bar{q}+1$, $\bar{q}^{3}$, and $\bar{q}^{3}+1$.

Let $f_{0}(X,Y,Z), f_{1}(X,Y,Z),\ldots,f_{v}(X,Y,Z)$ be distinct monomials of the same degree $m$, and such that $X$ is not a common factor of $f_{0}, f_{1},\ldots, f_{v}$.
Consider the set of rational functions
\begin{equation}\label{funz_raz}
\cL=\left\{\sum_{i=0}^{v}l_{i}f_{i}(x,y,z) : (l_{0}:l_{1}:\ldots:l_{v})\in\PP^{v}(\FF_{q^{2}})\right\}.
\end{equation}

Note that the second equation of $\phi(\cX)$
and $\sum_{i=0}^{v}l_{i}f_{i}(X,Y,Z)=0$
are the equations of two cones with the same vertex $O=(1:0:0:0)$. Hence, the zeros of any $\alpha\in\cL$ lie in the common lines of the two cones.
Let $\pi$ be the plane with equation $X=T$. Denote by $\cC_{1}$ and $\cC_{2}$ the plane curves obtained as the intersection of $\pi$ with the two cones. Clearly, the common lines of the two cones are the lines joining $O$ to the intersection points of $\cC_{1}$ and $\cC_{2}$. Assume that $\alpha=\sum_{i=0}^{v}l_{i}f_{i}(x,y,z)$. Then, in the $YZ$-plane the equations of $\cC_{1}$ and $\cC_{2}$ are
\begin{equation}\label{curvacono1}
\mathcal{C}_{1} : \left(a+b^{\bar{q}}Y-Z\right)^{\bar{q}}+(a+b^{\bar{q}}Y-Z)=(b+Y)^{\bar{q}+1}
\end{equation}
and
\begin{equation}\label{curvacono2}
\mathcal{C}_{2} : \displaystyle{\sum_{i=0}^{v}l_{i} f_{i}(1,Y,Z)=0}.
\end{equation}
Since these curves have order $\bar{q}+1$ and $m$, they have $m(\bar{q}+1)$ common points (not necessary distinct). Hence, the common lines of the cones are $m(\bar{q}+1)$. Note that the number of affine common points of $\cX$ and one of these lines $\l$ is either $\bar{q}^{2}-\bar{q}$ or $\bar{q}^{2}-\bar{q}+1$ according to whether $\l$ passes through $(1:1:0:0)$ or not.
Let $N$ be the number of affine zeros of the function $\alpha$.
Then,
$$N=(\bar{q}^{2}-\bar{q}+1)(m(\bar{q}+1)-M)+(\bar{q}^{2}-\bar{q})M,$$
where $M$ is the intersection multiplicity of $\cC_{1}$ and $\cC_{2}$ at the origin of the $YZ$-plane.

Consider now the morphism
$$\eta : \cC_{1}\rightarrow\PP^{v}(\bar{\FF}_{q^{2}}),\,\,\,\,\eta=(g_{0}:g_{1}:\ldots:g_{v}),$$
where $g_{i}=f_{i}(1,y,z)$. Note that if $g_{0},g_{1},\ldots, g_{v}$ are $\bar{\FF}_{q^{2}}$-linearly indipendent in the function field of $\cC_{1}$, then the morphism $\eta$ is non-degenerate. Let
$$\cD_{\eta}=\left\{E+div(\sum_{i=0}^{v}l_{i}g_{i}) : (l_{0}:l_{1}:\ldots: l_{v})\in\PP^{v}(\K)\right\}$$
be the linear series associated to $\eta$, where $E$ is the divisor of $\cC_{1}$ such that
$$v_{Q}(E)=-min\left\{v_{Q}(g_0), v_{Q}(g_1),\ldots, v_{Q}(g_v)\right\},$$
for any $Q\in\cC_{1}$. Therefore, we have $v+1$ distinct values for $M$, namely the integers $v_{P}(E)+j_{i}(P)$, $i=0,1,\ldots,v$, where $(j_{0}(P),j_{1}(P),\ldots,j_{v}(P))$ is the order sequence at $P$ of the morphism $\eta$.

Then the following result is obtained.

\begin{theorem}\label{W}
Let $P=(a,b,c)$ be an $\FF_{q^{2}}$-rational point of the GK-curve $\cX$, with $c\neq0$. Let $\cC_{1}$ be the plane curve with equation \eqref{curvacono1}, and let $g_{0},g_{1},\ldots,g_{v}$ be monomial functions in $y$, $z$, which are $\bar{\FF}_{q^{2}}$-linearly indipendent in the function field of $\cC_{1}$. Let $m$ be the maximum degree of $g_{0},g_{1},\ldots,g_{v}$, and let $v_{O}(E)=-min\left\{v_{O}(g_0), v_{O}(g_1),\ldots, v_{O}(g_v)\right\}$, where $O$ is the origin of the $YZ$-plane.
Then there exist $v+1$ non-gaps at $P$, say $N_{0}, N_{1},\ldots,N_{v}$, such that $m(\bar{q}^{3}-\bar{q})\leq N_{i}\leq m(\bar{q}^{3}+1)+v_{O}(E)$. More precisely $N_{i}=m(\bar{q}^{3}+1)-M_{i}$, where $M_{0},M_{1},\ldots,M_{v}$ are the intersection multiplicities at $O$ of $\cC_{1}$ and the plane curves with equation $\sum_{i=0}^{v}l_{i}g_{i}(Y,Z)=0$.
\end{theorem}

Thanks to Theorem \ref{W},  we are in a position  to obtain a description of the Weierstrass semigroup at an $\FF_{q^{2}}$-rational point of the GK-curve $\cX$ for $\bar{q}=2$ and $\bar{q}=3$.

\subsection{$\bar{q}=2$}\label{n2}

In this case, $\cX$ has affine equations
\begin{displaymath}
\left\{\begin{array}{ll}
Z^{3}=Y\left(1+X+X^{2}\right) \\
X^{2}+X=Y^{3}
\end{array}\right.,
\end{displaymath}
and $g=10$ is the genus of $\cX$.
Let $P=(a,b,c)$ be an $\FF_{q^{2}}$-rational point of $\cX$ such that $c\neq0$. Then,
the Weierstrass semigroup at $P$ coincides with the Weierstrass semigroup at $X_\infty$ of the curve $\phi(\cX)$ with equations
$$\displaystyle{\left(\frac{T}{c^{8}}+cX+\left(\frac{b^{2}-b^{8}}{c^{8}}\right)Y-\frac{Z}{c^{8}}\right)^{3}=}$$
$$\left(bX+Y\right) h(aX+b^{2}Y -Z)$$
 and
$$X\left(aX+b^{2}Y-Z\right)^{2}+X^{2}\left(aX+b^{2}Y-Z\right)=\left(bX+Y\right)^{3}.$$

By (\ref{ng09}),
$$div_{\infty}(x)=9X_\infty,\,\,\,\,div_{\infty}(y)=8X_\infty,\,\,\,\,div_{\infty}(z)=7X_\infty,$$
that is, $7$, $8$, and $9$ are non-gaps at $P$.
Let $\Upsilon$ be the semigroup generated by $7$, $8$, and $9$.
Since $\#(\Upsilon\cap\left[0,2g-1\right])=9<10=g$, $H(P)$ is larger than $\Upsilon$. Hence, there is precisely one non-gap not bigger than $2g-1=19$ that does not belong to $\Upsilon$.

Let
$$f_{0}=XZ,\,\,\,\,f_{1}=Z^{2},\,\,\,\,f_{2}=Y^{2},\,\,\,\,f_{3}=YZ,$$
and let $\eta : \cC_{1}\rightarrow\PP^{3}(\bar{\FF}_{q^{2}})$, $\eta=(z:z^{2}:y^{2}:yz)$; here, $\cC_{1}$ is the plane curve with equation $\cC_{1} : Z+Z^{2}+c^{3}Y^{2}+Y^{3}=0$. Since $z,z^{2},y^{2},yz$ are $\bar{\FF}_{q^{2}}$-linearly indipendent in the function field of $\cC_{1}$, the morphism $\eta$ is non-degenerate.
Consider now the linear series associated to $\eta$:
$$\cD_{\eta}=\left\{E+div(\sum_{i=0}^{i=3}l_{i}g_{i}) | (l_{0}:l_{1}:l_{2}:l_{3})\in\PP^{3}(\bar{\FF}_{q^{2}})\right\},$$
where $g_{i}(y,z)=f_{i}(1,y,z)$, and $E$ is a divisor of $\cC_{1}$ such that
$$v_{Q}(E)=-min\left\{v_{Q}(z), v_{Q}(z^{2}), v_{Q}(y^{2}), v_{Q}(yz)\right\},$$
for any $Q\in\cC_{1}$. Note that the line with equation $Z=0$ is the tangent line to $\cC_{1}$ at the origin $O$ of the $YZ$-plane, and its intersection multiplicity with $\cC_{1}$ is $2$. Hence, $v_{O}(E)=-2$ holds. Therefore, from Theorem~\ref{W} there exist four non-gaps
in $I=\left[12, 16\right]$.
Taking into account that $14,15,16$ are precisely the integers in $I\cap\left\langle 7,8,9\right\rangle$, then the non-gap $N$ is either $12$ or $13$.
According to Theorem \ref{W}, we implemented a standard intersection multiplicity algorithm in order to compute the intersection multiplicities at $O$ of $\cC_1$ and the plane curves with equation $\sum_{i=0}^{3}l_{i}g_{i}(Y,Z)=0$.
It turned out that $N=13$, and that a rational function $\beta$ such that $div_\infty(\beta)=13X_\infty$ was
$$\beta=xz+\left(\frac{1}{c^{9}}+1\right)z^{2}+c^{3}y^{2}+\frac{1}{c^{3}}yz.$$

It is easly seen that any integer greater than $20$ belongs to the numerical semigroup $\Delta$ generated by $7,8,9,13$. Hence, the genus of $\Delta$ equals $g$, and the following result is obtained.

\begin{theorem}\label{brrrr}
Let $P=(a,b,c)$ be an $\FF_{8^{2}}$-rational point of the GK-curve $\cX$, such that $c\neq0$. Then the Weierstrass semigroup of $\cX$ at $P$ is generated by $7,8,9$, and $13$. Moreover,
\begin{eqnarray*}
div_\infty( {\bar x})=9P,\qquad
div_\infty({\bar y})=8P,\qquad
div_\infty({\bar z})=7P,\\
div_\infty\left({\bar z}{\bar
x}+\left(\frac{1}{c^{9}}+1\right){\bar z}^{2}+c^{3}{\bar
y}^{2}+\frac{1}{c^{3}}{\bar y}{\bar z} \right)=13P
\end{eqnarray*}
where $\,\,\,{\bar x}=\frac{1}{-a^{8}-x+b^{8}y+c^{8}z},\,\,\,\,\,{\bar
y}=\frac{-b+y}{-a^{8}-x+b^{8}y+c^{8}z},\,\,\,\,$ and\\
${\bar z}=\frac{-a^{2}-x+b^{2}y}{-a^{8}-x+b^{8}y+c^{8}z}$.
\end{theorem}

\subsection{$\bar{q}=3$}

Let $\bar{q}=3$. In this case, $\cX$ has affine equations
\begin{displaymath}
\left\{\begin{array}{ll}
Z^{7}=Y\left(2+X^{2}+2X^{4}+X^{6}\right) \\
X^{3}+X=Y^{4}
\end{array}\right.,
\end{displaymath}
and $g=99$.
Let $P=(a,b,c)$ be an $\FF_{q^{2}}$-rational point of $\cX$ such that $c\neq0$.
Then,
the Weierstrass semigroup at $P$ coincides with the Weierstrass semigroup at $X_\infty$ of the curve $\phi(\cX)$ with equations
$$\displaystyle{\left(\frac{T}{c^{27}}+cX+\left(\frac{b^{3}-b^{27}}{c^{27}}\right)Y-\frac{Z}{c^{27}}\right)^{7} =}$$
$$\left(bX+Y\right) h(aX+b^{3}Y-Z)$$
and
$$X\left(aX+b^{3}Y-Z\right)^{3}+X^{3}\left(aX+b^{3}Y-Z\right)=\left(bX+Y\right)^{4}.$$

By (\ref{ng09}),
$$div_{\infty}(x)=28X_\infty,\,\,\,\,div_{\infty}(y)=27X_\infty,\,\,\,\,div_{\infty}(z)=25X_\infty,$$
that is $25$, $27$, and $28$ are non-gaps at $P$.
Let $\Upsilon$ be the semigroup generated by $25$, $27$, and $28$.
Since $\#(\Upsilon\cap\left[0,2g-1\right])=85<99=g$, $H(P)$ is larger than $\Upsilon$.

As for the case $\bar{q}=2$, we use Theorem \ref{W} to find the non-gaps at $P$ that do not belong to $\Upsilon$.
Let
$$f_0=X^{2}Z,\quad f_1=XZ^{2},\quad f_2=XYZ,$$
$$f_{3}=Y^{2}Z,\quad f_{4}=Z^{3},\quad f_{5}=YZ^{2},\quad f_{6}=Y^{3};$$
here, $\cC_{1} : 2Z+2Z^{3}+c^{7}Y^{3}+2Y^{4}=0$.
Since $z,z^{2},yz,y^{2}z,z^{3},yz^{2},y^{3}$ are $\bar{\FF}_{q^{2}}$-linearly indipendent in the function field of $\cC_{1}$, the morphism $\eta=(z:z^{2}:yz:y^{2}z:z^{3}:yz^{2}:y^{3}) : \cC_{1}\rightarrow\PP^{6}(\bar{\FF}_{q^{2}})$ is non-degenerate.
Consider now the linear series associated to $\eta$:
$$\cD_{\eta}=\left\{E+div(\sum_{i=0}^{i=6}l_{i}g_{i}) | (l_{0}:l_{1}:\ldots:l_{6})\in\PP^{6}(\bar{\FF}_{q^{2}})\right\},$$
where $g_{i}(y,z)=f_{i}(1,y,z)$, and $E$ is a divisor of $\cC_{1}$ such that
$$v_{Q}(E)=-min\left\{v_{Q}(g_{0}), v_{Q}(g_{1}),\ldots, v_{Q}(g_{6})\right\},$$
for any $Q\in\cC_{1}$. Note that the line with equation $Z=0$ is the tangent line to $\cC_{1}$ at the origin $O$ of the $YZ$-plane, and its intersection multiplicity with $\cC_{1}$ is $3$. Hence, $v_{O}(E)=-3$ holds. Therefore, by Theorem~\ref{W} there exist seven non-gaps
in $I=\left[72, 81\right]$.

As for the case $\bar{q}=2$, we implemented a standard intersection multiplicity algorithm in order to compute the intersection multiplicities at $O$ of $\cC_1$ and the plane curves with equation $\sum_{i=0}^{6}l_{i}g_{i}(Y,Z)=0$.
We obtained that $74$ was a non-gap at $X_\infty$  not belonging to
 $\Upsilon$, and that a rational function $\beta$ such that $div_\infty(\beta)=74X_\infty$ was
$$\displaystyle{\beta=x^2z+\frac{1}{c^{28}}xz^{2}+\frac{1}{c^{7}}xzy+\frac{1}{c^{14}}zy^{2}+}$$
$$\left(1+\frac{1}{c^{56}}\right)z^{3}+\frac{2}{c^{35}} z^{2}y+2c^{7}y^{3}.$$
By straightforward computation,  the number of integers less than $2g=198$ that belong to the subgroup generated by $25,27,28$, and $74$ is $96<g$. Therefore, some other non-gap is missing and we need to apply Theorem \ref{W} to another set of monomial functions.
Note that the following rational functions of $\cC_{1}$ are $\bar{\FF}_{q^{2}}$-linearly indipendent
$$y^{3},\,z,\,yz,\,y^{2}z,\,z^{2},\,y^{3}z,\,yz^{2},$$
$$y^{2}z^{2},\,z^{3},\,y^{3}z^{2},\,yz^{3},\,y^{2}z^{3},\,z^{4},\,yz^{4},\,z^{5}.$$
Consider the morphism associated to these rational functions.
Arguing as before, we applied a standard intersection algorithm to find
a non-gap  not belonging to the semigroup generated by
$25,27,28$, and $74$. A non-gap with this property turned out to be $121$. 
By straightforward computation it is easily seen that the genus of the numerical semigroup $\Delta$ generated by $25,27,28,74,$ and $121$ is $99$. Therefore, $\Delta$ coincides with the Weierstrass semigroup of $\varphi(\cX)$ at $X_\infty$.

An explicit description of a rational function $\gamma$ such that
$div_\infty(\gamma)=121X_\infty$ for a generic choice $(a,b,c)$
seems to be difficult to achieve. The intersection multiplicity algorithm  provided  such a function for a
specific choice of  $(a,b,c)\in \cO_2$ (note that this is not a
restriction, since by Theorem \ref{action} the automorphism group $Aut(\cX)$ acts transitively
on $\cO_2$). Let $\omega$ be an element of $\FF_{27^{2}}$ such that
$ \omega^6-\omega^4+\omega^2-\omega-1=0$. Then $(\omega^{11},
\omega^{280},\omega^{88})$ is a point in $\cO_2$, and
\begin{eqnarray*}
\gamma=\omega^{588}{  x}^2{  y}^3 + \omega^{336}{  x}^4{  z} +
\omega^{448}{  x}^3{  y}{  z} + \omega^{560}{
x}^2{  y}^2{  z} + \\
\omega^{700}{  x}^3{  z}^2 + 
\omega^{112}{  x}{  y}^3{  z}  + \omega^{112}{  x}^2{  y}{ z}^2  +
    \omega^{84}{  x}{  y}^2{  z}^2  + \\
    \omega^{196}{  x}^2{  z}^3  +
    2{  y}^3{  z}^2  +
    \omega^{392}{  x}{  y}{  z}^3  +
    \omega^{28}{  y}^2{  z}^3  + \\
    \omega^{504}{  x}{  z}^4
    + \omega^{644}{  y}{  z}^4  +
    \omega^{280}{  z}^5
\end{eqnarray*}
is such that $div_\infty(\gamma)=121X_\infty$.

Therefore, the following result is obtained.

\begin{theorem}\label{br2}
Let $\omega$ be an element of $\FF_{27^{2}}$ such that
$$
\omega^6-\omega^4+\omega^2-\omega-1=0.
$$
Let $P=(a,b,c)$ be an $\FF_{27^{2}}$-rational point of the
GK-curve $\cX$, such that $c\neq0$. Then, the Weierstrass
semigroup of $\cX$ at $P$ is generated by $25,27,28,74$, and
$121$. Moreover,
$$div_\infty( {\bar x})=28P,\quad
div_\infty({\bar y})=27P,$$
$$div_\infty({\bar z})=25P,\quad div_\infty\left(\bar{\beta}\right)=74P$$
where ${\bar x}=\frac{1}{-a^{27}-x+b^{27}y+c^{27}z}$, ${\bar
y}=\frac{-b+y}{-a^{27}-x+b^{27}y+c^{27}z}$, ${\bar
z}=\frac{-a^{3}-x+b^{3}y}{-a^{27}-x+b^{27}y+c^{27}z}$, and
$$\bar{\beta}={\bar z}{\bar x}^2+\frac{1}{c^{28}}{\bar x}{\bar z}^{2}+\frac{1}{c^{7}}{\bar x}{\bar y}{\bar z}+\frac{1}{c^{14}}{\bar z}{\bar y}^{2}+$$
$$\left(1+\frac{1}{c^{56}}\right){\bar z}^{3}+\frac{2}{c^{35}} {\bar z}^{2}{\bar y}+2c^{7}{\bar y}^{3}.$$
 When
$a=\omega^{11}$, $b=\omega^{280}$, $c=\omega^{88}$ then
$$
div_\infty(\bar\gamma)=121P
$$
where
\begin{eqnarray*}
\bar\gamma=\omega^{588}{\bar x}^2{\bar y}^3 + \omega^{336}{\bar
x}^4{\bar z} + \omega^{448}{\bar x}^3{\bar y}{\bar z} +
\omega^{560}{\bar x}^2{\bar y}^2{\bar z} + \\
\omega^{700}{\bar
x}^3{\bar z}^2 +  \omega^{112}{\bar x}{\bar y}^3{\bar z}  +
\omega^{112}{\bar x}^2{\bar y}{\bar z}^2  +
    \omega^{84}{\bar x}{\bar y}^2{\bar z}^2  + \\
    \omega^{196}{\bar x}^2{\bar z}^3  +
    2{\bar y}^3{\bar z}^2  +
    \omega^{392}{\bar x}{\bar y}{\bar z}^3  +
    \omega^{28}{\bar y}^2{\bar z}^3  + \\
    \omega^{504}{\bar x}{\bar z}^4
    + \omega^{644}{\bar y}{\bar z}^4  +
    \omega^{280}{\bar z}^5.
\end{eqnarray*}
\end{theorem}

\section{AG Codes and Improved AG Codes associated to the GK-curve}\label{sec5}

Throughout this section we keep the notation of the previous Sections. For $q=\bar{q}^3$, let $\cX$ be the GK curve defined by  (\ref{equazioni}). Let $P$ be an $\FF_{q^2}$-rational point of $\cX$, and let $D$ be the divisor consisting of the sum of all the remaining $\FF_{q^2}$-rational points of $\cX$.
Let $C_{\l}(P)$ be the dual of the AG code $C(D,\rho_{\l}P)$, with $\rho_{\l}\in H(P)$. Let ${\tilde C}_d(P)$ be the improved AG code, as defined in Section \ref{sec2}.

The aim of this section is to determine the parameters of both codes $C_{\l}(P)$ and ${\tilde C}_d(P)$, and to compare such parameters with those of the known codes.
We apply Theorems \ref{dist} and \ref{dist2}. Note that the bounds appearing in the statements of Theorems \ref{dist} and \ref{dist2} depend only on  Weierstrass semigroup $H(P)$. As $H(P)$ is invariant under the action of $Aut(\cX)$, we only consider one point per orbit under the action of $Aut(\cX)$. Henceforth, we assume that $P_i$ is a point of $\cO_i$, for $i=1,2$.

Note that by Proposition \ref{0}, Equations (\ref{eqqu1}) and (\ref{eqqu2}), and Theorems \ref{brrrr} and \ref{br2}, for each $P_i$ and for every $m\in H(P_i)$ we can construct a rational function $f$ such that $div_\infty{f}=mP_i$. Therefore, it is possible to construct a parity check matrix for all codes $C_{\l}(P_i)$ and ${\tilde C}_d(P_i)$

\begin{remark}It is well-known, see e.g. \cite{ST2}, that the permutation automorphism group of a code $C_\l(P)$ contains a subgroup isomorphic to the stabilizer of $P$ in $Aut(\cX)$, provided that the length of the code is larger than $2g+2$. Then by Corollary \ref{stabilizer} the code $C_\l(P_1)$ has an automorphism group of size $\bar{q}^3(\bar{q}^2-1)(\bar{q}^2-\bar{q}+1)$ , whereas   $C_\l(P_2)$ has a cyclic automorphism  group of size $\bar{q}^2-\bar{q}+1$.
\end{remark}

\subsection{$\bar{q}=2$}
Tables \ref{orbita1}-\ref{param2} describe the parameters of the codes $C_\ell(P_1)$, $C_\ell(P_2)$, ${\tilde C}_d(P_1)$, ${\tilde C}_d(P_2)$; the entries can be easily deduced from  Proposition \ref{0} and Theorem  \ref{brrrr}. In some cases the entries in MinT's tables \cite{MINT} are improved.

%
%
Other improvements can be obtained by using the following propagation rules.
\begin{proposition}[see Exercise 7 in \cite{TV}]\label{spoiling}${}$
\begin{itemize}
\item If there is a $q$-ary linear code of lenght $n$, dimension $k$ and minimum distance $d$, then for each non-negative integer $s< d$ there exists a $q$-ary linear code of length $n$, dimension $k$ and minimum distance $d-s$.

\item If there is a $q$-ary linear code of lenght $n$, dimension $k$ and minimum distance $d$, then for each non-negative integer $s< k$ there exists a $q$-ary linear code of length $n$, dimension $k-s$ and minimum distance $d$.
\item If there is a $q$-ary linear code of lenght $n$, dimension $k$ and minimum distance $d$, then for each non-negative integer $s<k$ there exists a $q$-ary linear code of length $n-s$, dimension $k-s$ and minimum distance $d$.
\item If there is a $q$-ary linear code of lenght $n$, dimension $k$ and minimum distance $d$, then for each non-negative integer $s<\min\{n-k-1,d\}$ there exists a $q$-ary linear code of length $n-s$, dimension $k$ and minimum distance $d-s$.
\end{itemize}
\end{proposition}
Therefore, the following result is obtained.
\begin{theorem}\label{new}
Linear codes over $\FF_{64}$ with parameters as in Table {\rm
\ref{chiave}} exist.
\end{theorem}

\subsection{$\bar{q}=3$} Table \ref{tretre} describes
some of the codes  ${\tilde C}_d(P_i)$, $d\le 2g$, $i=1,2$, over
the field $\FF_{3^6}$. The parameters of these codes
can be easily obtained taking into account Proposition \ref{0} and
Theorem \ref{br2}, together with Proposition \ref{spoiling}.

\begin{table*}
\renewcommand{\arraystretch}{1.3}
\caption{Improvements on \cite{MINT} - $q=64$} \label{chiave}
\centering
\begin{tabular}{|c|c|c|c||c|c|c|c||c|c|c|c|}
\hline
n & k & d & {\rm Ref.} & n & k & d & {\rm Ref.}& n & k & d &{\rm Ref.}\\
\hline
$224$ & $204$ & $13$ & ${C}_{19}(P_2)$, ${\tilde C}_{13}(P_2)$ &  $223$ & $203$ & $13$ & Prop. \ref{spoiling} &$222$ & $202$ & $13$ & Prop. \ref{spoiling} \\
$221$ & $201$ & $13$ & Prop. \ref{spoiling} &$220$ & $200$ & $13$ & Prop. \ref{spoiling} &$219$ & $199$ & $13$ & Prop. \ref{spoiling} \\
$218$ & $198$ & $13$ & Prop. \ref{spoiling} &$217$ & $197$ & $13$ & Prop. \ref{spoiling} &$216$ & $196$ & $13$ & Prop. \ref{spoiling}  \\
$215$ & $195$ & $13$ & Prop. \ref{spoiling} &$214$ & $194$ & $13$ & Prop. \ref{spoiling} &$213$ & $193$ & $13$ & Prop. \ref{spoiling}\\
$212$ & $192$ & $13$ & Prop. \ref{spoiling} &$211$ & $191$ & $13$ & Prop. \ref{spoiling} &$210$ & $190$ & $13$ & Prop. \ref{spoiling}\\
$209$ & $189$ & $13$ &Prop. \ref{spoiling} &$208$ & $188$ & $13$ &Prop. \ref{spoiling} &$207$ & $187$ & $13$ & Prop. \ref{spoiling}\\
$206$ & $186$ & $13$ &Prop. \ref{spoiling} &$205$ & $185$ & $13$ &Prop. \ref{spoiling} &$204$ & $184$ & $13$ & Prop. \ref{spoiling} \\
$203$ & $183$ & $13$ &Prop. \ref{spoiling} &$202$ & $182$ & $13$ &Prop. \ref{spoiling} &$201$ & $181$ & $13$ & Prop. \ref{spoiling} \\
$200$ & $180$ & $13$ &Prop. \ref{spoiling} &$224$ & $202$ & $14$ & {${C}_{21}(P_1)$,  ${\tilde C}_{14}(P_i)$} &$223$ & $201$ & $14$ & Prop. \ref{spoiling} \\
$222$ & $200$ & $14$ &Prop. \ref{spoiling} &$221$ & $199$ & $14$ &Prop. \ref{spoiling} &$220$ & $198$ & $14$ & Prop. \ref{spoiling} \\
$219$ & $197$ & $14$ &Prop. \ref{spoiling} &$218$ & $196$ & $14$ &Prop. \ref{spoiling} &$217$ & $195$ & $14$ & Prop. \ref{spoiling} \\
$216$ & $194$ & $14$ &Prop. \ref{spoiling} &$215$ & $193$ & $14$ &Prop. \ref{spoiling} &$214$ & $192$ & $14$ & Prop. \ref{spoiling} \\
$213$ & $191$ & $14$ &Prop. \ref{spoiling} &$212$ & $190$ & $14$ &Prop. \ref{spoiling} &$211$ & $189$ & $14$  & Prop. \ref{spoiling} \\
$210$ & $188$ & $14$ &Prop. \ref{spoiling} &$224$ & $201$ & $15$ &{${\tilde C}_{15}(P_1)$ } &$223$ & $200$ & $15$  & Prop. \ref{spoiling} \\
$222$ & $199$ & $15$ &Prop. \ref{spoiling} &$221$ & $198$ & $15$ &Prop. \ref{spoiling} &$220$ & $197$ & $15$ & Prop. \ref{spoiling} \\
$219$ & $196$ & $15$ &Prop. \ref{spoiling} &$218$ & $195$ & $15$ &Prop. \ref{spoiling} &$217$ & $194$ & $15$ & Prop. \ref{spoiling}\\
$216$ & $193$ & $15$ &Prop. \ref{spoiling} &$215$ & $192$ & $15$ &Prop. \ref{spoiling} &$214$ & $191$ & $15$ & Prop. \ref{spoiling} \\
$213$ & $190$ & $15$ &Prop. \ref{spoiling} &$212$ & $189$ & $15$ &Prop. \ref{spoiling} &$211$ & $188$ & $15$ & Prop. \ref{spoiling} \\
$210$ & $187$ & $15$ &Prop. \ref{spoiling} &$224$ & $196$ & $20$ &${C}_{27}(P_i)$, ${\tilde C}_{20}(P_i)$   &$223$ & $195$ & $20$ & Prop. \ref{spoiling} \\
$222$ & $194$ & $20$ &Prop. \ref{spoiling} &$221$ & $193$ & $20$ &Prop. \ref{spoiling} &$220$ & $192$ & $20$ & Prop. \ref{spoiling}\\
$219$ & $191$ & $20$ &Prop. \ref{spoiling} &$218$ & $190$ & $20$ &Prop. \ref{spoiling} &$217$ & $189$ & $20$ & Prop. \ref{spoiling}\\
$216$ & $188$ & $20$ &Prop. \ref{spoiling} &$215$ & $187$ & $20$ &Prop. \ref{spoiling} &$214$ & $186$ & $20$ & Prop. \ref{spoiling} \\
$213$ & $185$ & $20$ &Prop. \ref{spoiling} &$212$ & $184$ & $20$ &Prop. \ref{spoiling} &$211$ & $183$ & $20$ & Prop. \ref{spoiling}  \\
$210$ & $182$ & $20$ &Prop. \ref{spoiling} & &  &  &{} & &  &{} & \\
\hline
\end{tabular}
\end{table*}

\begin{table}
\renewcommand{\arraystretch}{1.3}
\caption{Codes $C_{\l}(P_1)$ - $q=64$} \label{orbita1}
\centering
\begin{tabular}{|c|c|c|c|c|}
\hline
$n$ & $k$ & $\rho_{\ell}$ & $\nu_{\ell}$ & $d_{ORD}$\\
\hline
224 & 223 & 0 & 2 & 2 \\
224 & 222 &6 & 2 & 2 \\
224 & 222 &8 & 2 & 2 \\
224 & 220 &9 & 3 & 3 \\
224 & 219 &12 & 4 & 3 \\
224 & 218 &14 & 4 & 3 \\
224 & 217 &15 & 3 & 3 \\
224 & 216 &16 & 4 & 4 \\
224 & 215 &17 & 5 & 5 \\
224 & 214 &18 & 6 & 6 \\
224 & 213 &20 & 6 & 6 \\
224 & 212 &21 & 6 & 6 \\
224 & 211 &22 & 8 & 6 \\
224 & 210 &23 & 9 & 6 \\
224 & 209 &24 & 6 & 6 \\
224 & 208 &25 & 10 & 8\\
224 & 207 &26 & 8 & 8 \\
224 & 206 &27 & 9 & 9 \\
224 & 205 &28 & 12 & 12 \\
224 & 204 &29 & 13 & 12 \\
224 & 203 &30 & 12 & 12 \\
224 & 202 &31 & 15 & 14 \\
224 & 201 &32 & 14 & 14 \\
224 & 200 &33 & 15 & 15 \\
224 & 199 &34 & 16 & 16 \\
224 & 198 &35 & 17 & 17 \\
224 & 197 &36 & 18 & 18 \\
224 & 196 &37 & 20 & 20 \\
224 & 195 &38 & 20 & 20 \\
\hline
\end{tabular}
\end{table}

\begin{table}
\renewcommand{\arraystretch}{1.3}
\caption{Codes $C_{\l}(P_2)$ - $q=64$} \label{orbita2}
\centering
\begin{tabular}{|c|c|c|c|c|}
\hline
$n$ & $k$ & $\rho_{\ell}$ & $\nu_{\ell}$ & $d_{ORD}$\\
\hline
224 & 223 &0 & 2 & 2\\
224 & 222 &7 & 2 & 2\\
224 & 221 &8 & 2 & 2\\
224 & 220 &9 & 2 & 2\\
224 & 219 &13 & 3 & 3\\
224 & 218 &14 & 4 & 3\\
224 & 217 &15 & 5 & 3\\
224 & 216 &16 & 4 & 3\\
224 & 215 &17 & 3 & 3\\
224 & 214 &18 & 4 & 4\\
224 & 213 &20 & 6 & 6\\
224 & 212 &21 & 8 & 7\\
224 & 211 &22 & 8 & 7\\
224 & 210 &23 & 8 & 7\\
224 & 209 &24 & 8 & 7\\
224 & 208 &25 & 7 & 7\\
224 & 207 &26 & 8 & 8\\
224 & 206 &27 & 9 & 9\\
224 & 205 &28 & 12 & 12\\
224 & 204 &29 & 13 & 13\\
224 & 203 &30 & 14 & 13\\
224 & 202 &31 & 13 & 13\\
224 & 201 &32 & 14 & 14\\
224 & 200 &33 & 15 & 15\\
224 & 199 &34 & 16 & 16\\
224 & 198 &35 & 17 & 17\\
224 & 197 &36 & 18 & 18\\
224 & 196 &37 & 20 & 20\\
224 & 195 &38 & 20 & 20\\
\hline
\end{tabular}
\end{table}

\begin{table}
\renewcommand{\arraystretch}{1.3}
\caption{Codes ${\tilde C}_d(P_1)$ -  $q=64$} \label{param1}
\centering
\begin{tabular}{|c|c|c|c|}
\hline
$n$ & $d$ & $r_d$ & $k\geq$ \\
\hline
224 & 3 & 4 & 220\\
224 & 4 & 6 & 218\\
224 & 5 & 9 & 215\\
224 & 6 & 10 & 214\\
224 & 7 & 14 & 210\\
224 & 8 & 14 & 210\\
224 & 9 & 16 & 208\\
224 & 10 & 18 & 206\\
224 & 11 & 19 & 205\\
224 & 12 & 19 & 205\\
224 & 13& 21 & 203\\
224 & 14 & 22 &
202\\
224 & 15 & 23 & 201\\
224 & 16 & 25 & 199\\
224 & 17 & 26 & 198\\
224 & 18 & 27 & 197\\
224 & 19 & 28 & 196\\
224 & 20 & 28 & 196\\
\hline
\end{tabular}
\end{table}

\begin{table}
\renewcommand{\arraystretch}{1.3}
\caption{Codes ${\tilde C}_d(P_2)$ - $q=64$} \label{param2}
\centering
\begin{tabular}{|c|c|c|c|}
\hline
$n$ & $d$ & $r_d$ & $k\geq$  \\
\hline
224 & 3 & 5 & 219\\
224 & 4 & 7 & 217\\
224 & 5 & 10 & 214\\
224 & 6 & 11 & 213\\
224 & 7 & 12 & 212\\
224 & 8 & 13 & 211\\
224 & 9 & 18 & 206\\
224 & 10 & 19 & 205\\
224 & 11 & 19 & 205\\
224 & 12 & 19 & 205\\
224 & 13 & 20 & 204\\
224 & 14 & 22 &
202\\
224 & 15 & 24 & 200\\
224 & 16 & 25 & 199\\
224 & 17& 26 & 198\\
224 & 18 & 27 & 197\\
224 & 19 & 28 & 196\\
224 & 20 & 28 & 196\\
\hline
\end{tabular}
\end{table}

\begin{tiny}
\begin{table*}
\renewcommand{\arraystretch}{1.3}
\caption{$q=3^6$, $n=6075$} \label{tretre}
\centering
\begin{tabular}{|c|c|c||c|c|c||c|c|c||c|c|c||c|c|c|}
\hline $k$ & $d$ & Ref. & $k$ & $d$ & Ref. & $k$ & $d$ & Ref. &
$k$ & $d$ & Ref. \\ \hline 
$6074$ & 2 & ${\tilde C}_2(P_1)$ &
$6071$ & 3 & ${\tilde C}_3(P_1)$ & 
$6068$ & 4 & ${\tilde C}_4(P_1)$ & 
$6063$ & 5 & ${\tilde C}_5(P_1)$ \\
\hline 
$6062$ & 6 & ${\tilde C}_6(P_1)$ & 
$6055$ & 7 & ${\tilde C}_7(P_1)$ & 
$6053$ & 8 & ${\tilde C}_8(P_1)$ & 
$6048$ & 9 & ${\tilde C}_9(P_1)$ \\
\hline 
$6045$ & 10 & ${\tilde C}_{10}(P_1)$ & 
$6042$ & 11 & ${\tilde C}_{11}(P_1)$ & 
$6041$ & 12 & ${\tilde C}_{12}(P_1)$ &
$6032$ & 13 & ${\tilde C}_{13}(P_1)$ \\
\hline 
$6031$ & 14 & ${\tilde C}_{14}(P_1)$ & 
$6027$ & 15 & ${\tilde C}_{15}(P_1)$ & 
$6024$ & 16 & ${\tilde C}_{16}(P_1)$ &
$6020$ & 17 & ${\tilde C}_{17}(P_1)$ \\
\hline 
$6019$ & 18 & ${\tilde C}_{18}(P_1)$ & 
$6013$ & 19 & ${\tilde C}_{19}(P_1)$ & 
$6012$ & 20 & ${\tilde C}_{20}(P_1)$ &
$6008$ & 21 & ${\tilde C}_{21}(P_1)$ \\
\hline 
$6004$ & 22 & ${\tilde C}_{22}(P_1)$ & 
$6003$ & 23 & ${\tilde C}_{23}(P_1)$ & 
$6002$ & 24 & ${\tilde C}_{24}(P_1)$ &
$5996$ & 25 & ${\tilde C}_{25}(P_2)$ \\
\hline 
$5995$ & 26 & ${\tilde C}_{26}(P_1)$ & 
$5994$ & 27 & ${\tilde C}_{27}(P_1)$ & 
$5992$ & 28 & ${\tilde C}_{28}(P_1)$ &
$5987$ & 30 & ${\tilde C}_{30}(P_1)$ \\
\hline 
$5983$ & 32 & ${\tilde C}_{32}(P_1)$ & 
$5981$ & 33 & ${\tilde C}_{33}(P_1)$ & 
$5980$ & 34 & ${\tilde C}_{34}(P_1)$ &
$5979$ & 35 & ${\tilde C}_{35}(P_1)$ \\
\hline 
$5978$ & 36 & ${\tilde C}_{36}(P_1)$ & 
$5973$ & 37 & ${\tilde C}_{37}(P_1)$ & 
$5972$ & 38 & ${\tilde C}_{38}(P_1)$ &
$5970$ & 39 & ${\tilde C}_{39}(P_1)$ \\
\hline 
$5969$ & 40 & ${\tilde C}_{40}(P_1)$ & 
$5966$ & 42 & ${\tilde C}_{42}(P_1)$ & 
$5961$ & 44 & ${\tilde C}_{44}(P_1)$ &
$5960$ & 45 & ${\tilde C}_{45}(P_1)$ \\
\hline 
$5958$ & 46 & ${\tilde C}_{46}(P_1)$ & 
$5956$ & 48 & ${\tilde C}_{48}(P_1)$ & 
$5952$ & 50 & ${\tilde C}_{50}(P_2)$ &
$5951$ & 51 & ${\tilde C}_{51}(P_2)$ \\
\hline 
$5949$ & 52 & ${\tilde C}_{52}(P_2)$ & 
$5946$ & 53 & ${\tilde C}_{53}(P_1)$ & 
$5945$ & 54 & ${\tilde C}_{54}(P_1)$ &
$5942$ & 55 & ${\tilde C}_{55}(P_1)$ \\
\hline
 $5940$ & 56 & ${\tilde C}_{56}(P_1)$ & 
$5938$ & 57 & ${\tilde C}_{57}(P_1)$ & 
$5937$ & 60 & ${\tilde C}_{60}(P_1)$ &
$5932$ & 62 & ${\tilde C}_{62}(P_1)$ \\
\hline 
$5931$ & 63 & ${\tilde C}_{63}(P_1)$ & 
$5929$ & 64 & ${\tilde C}_{64}(P_1)$ & 
$5928$ & 65 & ${\tilde C}_{65}(P_1)$ &
$5927$ & 66 & ${\tilde C}_{66}(P_1)$ \\
\hline
 $5926$ & 68 & ${\tilde C}_{68}(P_1)$ & 
$5924$ & 69 & ${\tilde C}_{69}(P_1)$ & 
$5922$ & 70 & ${\tilde C}_{70}(P_1)$ &
$5919$ & 71 & ${\tilde C}_{71}(P_1)$ \\
\hline 
$5918$ & 72 & ${\tilde C}_{72}(P_1)$ & 
$5917$ & 74 & ${\tilde C}_{74}(P_2)$ & 
$5916$ & 75 & ${\tilde C}_{75}(P_2)$ &
$5915$ & 76 & ${\tilde C}_{76}(P_2)$ \\
\hline 
$5914$ & 77 & ${\tilde C}_{77}(P_2)$ & 
$5913$ & 78 & ${\tilde C}_{78}(P_2)$ & 
$5910$ & 79 & ${\tilde C}_{79}(P_2)$ &
$5908$ & 80 & ${\tilde C}_{80}(P_1)$ \\
\hline 
$5906$ & 81 & ${\tilde C}_{81}(P_1)$ & 
$5905$ & 82 & ${\tilde C}_{82}(P_1)$ & 
$5904$ & 83 & ${\tilde C}_{83}(P_1)$ &
$5902$ & 84 & ${\tilde C}_{84}(P_1)$ \\
\hline 
$5899$ & 85 & ${\tilde C}_{85}(P_1)$ & 
$5898$ & 86 & ${\tilde C}_{86}(P_1)$ & 
$5897$ & 90 & ${\tilde C}_{90}(P_1)$ &
$5894$ & 91 & ${\tilde C}_{91}(P_1)$ \\
\hline 
$5892$ & 92 & ${\tilde C}_{92}(P_1)$ & 
$5891$ & 94 & ${\tilde C}_{94}(P_1)$ & 
$5890$ & 96 & ${\tilde C}_{96}(P_2)$ &
$5889$ & 99 & ${\tilde C}_{99}(P_2)$ \\
\hline 
$5888$ & 100 & ${\tilde C}_{100}(P_2)$ & 
$5885$ & 101 & ${\tilde C}_{101}(P_2)$ & 
$5884$ & 102 & ${\tilde C}_{102}(P_2)$ &
$5880$ & 103 & ${\tilde C}_{103}(P_1)$ \\
\hline 
$5878$ & 104 & ${\tilde C}_{104}(P_1)$ & 
$5877$ & 105 & ${\tilde C}_{105}(P_1)$ & 
$5875$ & 106 & ${\tilde C}_{106}(P_1)$ &
$5874$ & 108 & ${\tilde C}_{108}(P_1)$ \\
\hline 
$5872$ & 109 & ${\tilde C}_{109}(P_1)$ & 
$5871$ & 110 & ${\tilde C}_{110}(P_1)$ & 
$5869$ & 111 & ${\tilde C}_{111}(P_1)$ &
$5868$ & 112 & ${\tilde C}_{112}(P_1)$ \\
\hline
 $5866$ & 114 & ${\tilde C}_{114}(P_1)$ & 
$5865$ & 115 & ${\tilde C}_{115}(P_1)$ & 
$5864$ & 117 & ${\tilde C}_{117}(P_1)$ &
$5863$ & 120 & ${\tilde C}_{120}(P_2)$ \\
\hline
 $5862$ & 121 & ${\tilde C}_{121}(P_2)$ & 
$5860$ & 124 & ${\tilde C}_{124}(P_2)$ &
 $5857$ & 125 & ${\tilde C}_{125}(P_2)$ &
$5854$ & 126 & ${\tilde C}_{126}(P_1)$ \\
\hline
 $5852$ & 128 & ${\tilde C}_{128}(P_1)$ & 
$5851$ & 129 & ${\tilde C}_{129}(P_1)$ &
 $5849$ & 130 & ${\tilde C}_{130}(P_1)$ &
$5848$ & 131 & ${\tilde C}_{131}(P_1)$ \\
\hline
 $5847$ & 132 & ${\tilde C}_{132}(P_1)$ & 
$5846$ & 133 & ${\tilde C}_{133}(P_1)$ & 
$5844$ & 134 & ${\tilde C}_{134}(P_1)$ &
$5843$ & 135 & ${\tilde C}_{135}(P_1)$ \\
\hline 
$5842$ & 136 & ${\tilde C}_{136}(P_1)$ & 
$5841$ & 137 & ${\tilde C}_{137}(P_1)$ & 
$5840$ & 138 & ${\tilde C}_{138}(P_1)$ &
$5838$ & 139 & ${\tilde C}_{139}(P_1)$ \\
\hline 
$5837$ & 140 & ${\tilde C}_{140}(P_1)$ & 
$5836$ & 144 & ${\tilde C}_{144}(P_1)$ & 
$5835$ & 146 & ${\tilde C}_{146}(P_2)$ &
$5832$ & 148 & ${\tilde C}_{148}(P_2)$ \\
\hline 
$5830$ & 149 & ${\tilde C}_{149}(P_2)$ & 
$5829$ & 150 & ${\tilde C}_{150}(P_1)$ & 
$5828$ & 151 & ${\tilde C}_{151}(P_1)$ &
$5827$ & 152 & ${\tilde C}_{152}(P_1)$ \\
\hline 
$5825$ & 153 & ${\tilde C}_{153}(P_1)$ & 
$5823$ & 154 & ${\tilde C}_{154}(P_1)$ & 
$5822$ & 156 & ${\tilde C}_{156}(P_1)$ &
$5821$ & 157 & ${\tilde C}_{157}(P_1)$ \\
\hline 
$5820$ & 158 & ${\tilde C}_{158}(P_1)$ & 
$5818$ & 159 & ${\tilde C}_{159}(P_1)$ & 
$5817$ & 160 & ${\tilde C}_{160}(P_1)$ &
$5816$ & 161 & ${\tilde C}_{161}(P_1)$ \\
\hline 
$5815$ & 162 & ${\tilde C}_{162}(P_1)$ & 
$5814$ & 163 & ${\tilde C}_{163}(P_1)$ & 
$5813$ & 164 & ${\tilde C}_{164}(P_1)$ &
$5812$ & 165 & ${\tilde C}_{165}(P_1)$ \\
\hline 
$5811$ & 166 & ${\tilde C}_{166}(P_1)$ & 
$5810$ & 167 & ${\tilde C}_{167}(P_1)$ & 
$5809$ & 168 & ${\tilde C}_{168}(P_1)$ &
$5808$ & 171 & ${\tilde C}_{171}(P_1)$ \\
\hline
$5806$ & 172 & ${\tilde C}_{172}(P_1)$ & 
$5805$ & 173 & ${\tilde C}_{173}(P_2)$ & 
$5804$ & 174 & ${\tilde C}_{174}(P_2)$ &
$5802$ & 175 & ${\tilde C}_{175}(P_1)$ \\
\hline 
$5801$ & 177 & ${\tilde C}_{177}(P_1)$ & 
$5800$ & 178 & ${\tilde C}_{178}(P_1)$ & 
$5798$ & 179 & ${\tilde C}_{179}(P_1)$ &
$5797$ & 180 & ${\tilde C}_{180}(P_1)$ \\
\hline 
$5796$ & 181 & ${\tilde C}_{181}(P_1)$ & 
$5795$ & 182 & ${\tilde C}_{182}(P_1)$ & 
$5794$ & 183 & ${\tilde C}_{183}(P_1)$ &
$5793$ & 184 & ${\tilde C}_{184}(P_1)$ \\
\hline 
$5792$ & 185 & ${\tilde C}_{185}(P_1)$ & 
$5791$ & 186 & ${\tilde C}_{186}(P_1)$ & 
$5790$ & 187 & ${\tilde C}_{187}(P_1)$ &
$5789$ & 188 & ${\tilde C}_{188}(P_1)$ \\
\hline 
$5788$ & 189 & ${\tilde C}_{189}(P_1)$ & 
$5787$ & 190 & ${\tilde C}_{190}(P_1)$ & 
$5786$ & 191 & ${\tilde C}_{191}(P_1)$ &
$5785$ & 192 & ${\tilde C}_{192}(P_1)$ \\
\hline 
$5784$ & 193 & ${\tilde C}_{193}(P_1)$ & 
$5783$ & 194 & ${\tilde C}_{194}(P_1)$ & 
$5782$ & 195 & ${\tilde C}_{195}(P_1)$ &
$5781$ & 196 & ${\tilde C}_{196}(P_1)$ \\
\hline 
$5780$ & 198 & ${\tilde C}_{198}(P_1)$ &  &  &
 &  &  &  &
 &  &  \\
\hline
\end{tabular}
\end{table*}
\end{tiny}


\begin{thebibliography}{1}

\bibitem{MINT} MinT. (2009, January). Tables of optimal parameters for linear codes, 
University of Salzburg. Available: {http://mint.sbg.ac.at/}.


\bibitem{GO1} V.D. Goppa,  ``Codes on algebraic curves,'' {\em Dokl. Akad. NAUK, SSSR}, vol. {259}, pp. 1289--1290, 1981.

\bibitem{GO2} V.D. Goppa, ``Algebraic-geometric Codes,'' {\em Izv. Akad. NAUK, SSSR}, vol. 46, pp. 75--91,  1982.

\bibitem{HA} J.P. Hansen, ``Codes on the Klein Quartic, Ideals and decoding,'' {\em IEEE Trans. Inf. Theory}, vol. 33, pp. 923--925, 1987.

\bibitem{HS} J.P. Hansen and H. Stichtenoth, ``Group Codes on Certain Algebraic Curves with Many Rational Points,'' {\em AAECC}, vol. 1, pp. 67--77, 1990.

\bibitem{M0} G.L. Matthews, ``Codes from the Suzuki function field,'' {\em IEEE Trans. Inf. Theory}, vol. 50, pp. 3298--3302, 2004.

\bibitem{M1} G.L. Matthews, ``Weierstrass Semigroups and Codes from a Quotient of the Hermitian Curve,'' {\em Des. Codes Cryptogr.}, vol.  37, pp.  473--492, 2005.


\bibitem{ST} H. Stichtenoth,  ``A note on Hermitian Codes over $GF(q^2)$,'' {\em IEEE Trans. Inf. Theory}, vol. 34, pp. 1345--1348, 1988.

\bibitem{TI} H.J. Tiersma,  ``Remarks on codes from Hermitian curves,'' {\em IEEE Trans. Inf. Theory}, vol.  33, pp. 605--609, 1987.

\bibitem{XC} C.P. Xing and H. Chen, ``Improvements on parameters of one-point AG codes from Hermitian curves,'' {\em IEEE Trans. Inf. Theory}, vol. 48, pp. 535--537, 2002.

\bibitem{XL} C.P. Xing and S. Ling, ``A class of Linear Codes with Good Parameters from Algebraic Curves,''  {\em IEEE Trans. Inf. Theory}, vol. 46, pp. 1527--1532, 2000.


\bibitem{YK} K. Yang and P.V. Kumar, ``On the true minimum distance of Hermitian codes,'' {\em Coding theory and
algebraic geometry, Lecture Notes in Math.},  vol. 1518, pp. 99-107, 1992.


\bibitem{GK} M. Giulietti and G. Korchm\'aros, ``A new family of maximal curves over a finite field,'' \emph{Mathematische Annalen}, vol.  343, pp. 229--245, 2009.

\bibitem{HLP} T. H{\o}holdt, J.H. van Lint and R. Pellikaan, ``Algebraic Geometry codes,'' in \emph{Handbook of Coding Theory}, vol. 1,  V.S. Pless, W.C.Huffman and R.A. Brualdi Eds. Amsterdam: Elsevier, 1998, pp. 871--961.


\bibitem{GKL} A. Garcia, S.J. Kim and R.F. Lax, ``Consecutive Weierstrass gaps and minimum distance of Goppa codes,'' \emph{J. Pure Appl. Algebra}, vol. 84, pp. 199--207, 1993.

\bibitem{SV} K.O. St\"ohr and J.F. Voloch,  ``Weierstrass points and curves over finite fields,'' \emph{Proc. London Math. Soc. $(3)$}, vol. 52, pp. 1--19, 1986.





\bibitem{FGT} R. Fuhrmann, A. Garcia and F. Torres,  ``On maximal curves,'' \emph{J. Number Theory}, vol. 67, pp. 29--51, 1997.




\bibitem{HKT}  J.W.P. Hirschfeld, G. Korchm\'aros and F. Torres, {\em Algebraic Curves over a Finite Field}, Princeton/Oxford: Princeton University Press, 2008.



\bibitem{ST2} H. Stichtenoth, ``On automorphisms of geometric Goppa codes,'' {\em J. Algebra}, vol. 130, pp. 113-121, 1990.





\bibitem{TV} M.A. Tsfasman and S.G. Vladut, {\em Algebraic-goemetric codes}, Amsterdam: Kluwer,  1991.



\end{thebibliography}
\end{document}